\def\vers{May 31, 2004, v.4}
\documentstyle{amsppt}
\magnification=1200
\hsize=6.5truein
\vsize=8.9truein
\hoffset=-10pt
\voffset=-20pt
\topmatter
\title Some Consequences of Perversity
of Vanishing Cycles
\endtitle
\author Alexandru Dimca and Morihiko Saito
\endauthor
\abstract
For a holomorphic function on a complex manifold, we show
that the vanishing cohomology of lower degree at a point
is determined by that for the points near it, using the
perversity of the vanishing cycle complex.
We calculate it explicitly in the case the hypersurface
has simple normal crossings outside the point.
We also give some applications to the monodromy.
\endabstract
\endtopmatter
\tolerance=1000
\baselineskip=12pt
\document
\def\scirc{\raise.2ex\hbox{${\scriptstyle\circ}$}}
\def\ssbull{\raise.2ex\hbox{${\scriptscriptstyle\bullet}$}}
\def\msum{\hbox{$\sum$}}
\def\mprod{\hbox{$\prod$}}
\def\mopls{\hbox{$\bigoplus$}}
\def\motim{\hbox{$\bigotimes$}}
\def\mcap{\hbox{$\bigcap$}}
\def\mcup{\hbox{$\bigcup$}}
\def\bC{{\Bbb C}}
\def\bD{{\Bbb D}}
\def\bP{{\Bbb P}}
\def\bQ{{\Bbb Q}}
\def\boH{\bold{H}}
\def\boR{\bold{R}}

\def\cF{{\Cal F}}
\def\cH{{\Cal H}}
\def\cO{{\Cal O}}
\def\tH{{\widetilde{H}}}
\def\tK{{\widetilde{K}}}

\def\oY{{\overline{Y}}}
\def\olam{{\overline{\lambda}}}
\def\oSig{{\overline{\Sigma}}}
\def\red{\text{{\rm red}}}
\def\ver{\text{{\rm ver}}}
\def\hor{\text{{\rm hor}}}
\def\Gr{\hbox{{\rm Gr}}}
\def\Ker{\hbox{{\rm Ker}}}
\def\Im{\hbox{{\rm Im}}}
\def\IC{\hbox{{\rm IC}}}
\def\div{\hbox{{\rm div}}\,}
\def\can{\hbox{{\rm can}}}
\def\Var{\hbox{{\rm Var}}}

\def\supp{\hbox{{\rm supp}}\,}
\def\Sing{\hbox{{\rm Sing}}\,}

\def\simto{\buildrel\sim\over\longrightarrow}
\def\SameAuthor{\vrule height3pt depth-2.5pt width1cm}

\centerline{\bf Introduction}\footnote""{{\it Date}\,: \vers}

\bigskip\noindent
Let
$ f $ be a nonconstant holomorphic function on a complex
analytic space
$ X $.
For each
$ x \in Y := f^{-1}(0) $,
we have the vanishing cohomology
$ \tH^{j}(F_{x},\bQ) $ where
$ F_{x} $ denotes the (typical) fiber of the Milnor
fibration around
$ x $, and
$ \tH $ means the reduced cohomology.
It has been observed by many people that there are certain
relations between the
$ \tH^{j}(F_{x},\bQ) $ for
$ x \in Y $.
It is well-known that they form a constructible sheaf on
$ Y $ (called the vanishing cohomology sheaf).
P. Deligne [7] constructed a sheaf complex
$ \varphi_{f}\bQ_{X} $ on
$ Y $ (called the vanishing cycle complex) such that its
cohomology sheaves
$ \cH^{j}\varphi_{f}\bQ_{X} $ are the vanishing
cohomology sheaves.

Let
$ L_{x} $ denote the intersection of
$ Y $ with a sufficiently small
sphere around
$ x \in Y $ (in a smooth ambient space), which is called the
{\it link} of
$ \{x\} $ in
$ Y $.
Let
$ T_{u} $,
$ T_{s} $ be respectively the unipotent and semisimple part
of the monodromy
$ T $, and put
$ N = \log T_{u} $.
Let
$ \tH^{n-1}(F_{x},\bQ)_{1} $ and
$ \tH^{n-1}(F_{x},\bQ)_{\ne 1} $ denote the unipotent
and non unipotent monodromy part, which are defined by
$ \Ker\,(T_{s}-1) $ and
$ \mopls_{\lambda\ne 1}\Ker\,(T_{s}-\lambda) $ (after a
scalar extension) respectively,
and similarly for the cohomology with compact supports.

\medskip\noindent
{\bf 0.1.~Theorem.}
{\it Assume that
$ \bQ_{X}[n+1] $ is a perverse sheaf {\rm (}e.g.
$ X $ is a locally complete intersection of dimension
$ n + 1) $, and
$ n \ge 1 $.
Then there are canonical isomorphisms
$$
\tH^{j}(F_{x},\bQ) = \boH^{j}(L_{x},
\varphi_{f}\bQ_{X}|_{L_{x}})\quad\text{for}\,\,j < n - 1,
$$
and a short exact sequence
$$
0 \to \tH^{n-1}(F_{x},\bQ) \to \boH^{n-1}
(L_{x},\varphi_{f}\bQ_{X}|_{L_{x}}) \to K_{x} \to 0.
$$
Here
$ K_{x} $ is the kernel of a morphism
$ \beta_{\varphi} $ which is the direct sum of
$$
\aligned
\beta_{\varphi,1} :
H_{c}^{n}(F_{x},\bQ)_{1}(-1)
&\to
H^{n}(F_{x},\bQ)_{1},
\\
\beta_{\varphi,\ne 1} :
H_{c}^{n}(F_{x},\bQ)_{\ne 1}
&\to
H^{n}(F_{x},\bQ)_{\ne 1},
\endaligned
$$
where
$ (-1) $ denotes the Tate twist, and
$ \beta_{\varphi,\ne 1} $ coincides with the natural
morphism {\rm (}i.e. corresponds to the natural
intersection form if
$ X $ is a rational homology manifold\,{\rm ).}
If
$ X $ is a rational homology manifold at
$ x $, then
$ N\beta_{\varphi,1} $ coincides with the natural morphism.
These morphisms and the short exact sequence are compatible
with mixed Hodge structure.
}

\medskip
In the
$ 1 $-dimensional singular locus case, a similar assertion
was obtained in [18], [19], see also [1].
Theorem (0.1) means that
$ \tH^{j}(F_{x},\bQ) $ for
$ j < n - 1 $ (resp.
$ j = n - 1) $ is completely (resp.
partially) determined by the restriction of
$ \varphi_{f}\bQ_{X} $ to the complement of
$ x $, and only
$ \tH^{n}(F_{x},\bQ) $ is essentially
interesting if we know well about the restriction of
$ \varphi_{f}\bQ_{X} $ to the complement of
$ x $.
The proof easily follows from the well-known fact that the
vanishing cycle complex
$ \varphi_{f}\bQ_{X} $ is a (shifted) perverse sheaf.
Actually, the first two assertions of Theorem (0.1) are
essentially equivalent to the perversity of
$ \varphi_{f}\bQ_{X} $, assuming the perversity of its
restriction to the complement of
$ x $.
The hypercohomology
$ \boH^{j}(L_{x},\varphi_{f}\bQ_{X}|_{L_{x}})
$ can be calculated by using spectral sequences
(2.2--3).
The mixed Hodge structure on
$ H^{j}(F_{x},\bQ) $ can be calculated by using the
weight spectral sequence (1.5), see also [14] for the
unipotent monodromy case, and [20] for the isolated
singularity case.

In Theorem (0.1) we can replace the vanishing cycle complex
$ \varphi_{f}\bQ_{X} $ with the nearby cycle complex
$ \psi_{f}\bQ_{X} $ in [7], and
$ \beta_{\varphi} $ with
$ \beta_{\psi} : H_{c}^{n}(F_{x},\bQ) \to H^{n}(F_{x},\bQ) $.
In this case
$ \beta_{\psi} $
is a natural morphism, and in the isolated
singularity case (where
$ X $ is smooth), we get a well-known relation between
the cohomology of the Milnor fiber and the link.
Note that the morphism
$ \beta_{\varphi} $ in Theorem (0.1) for
$ \varphi $ in the isolated singularity case is an
isomorphism (i.e. the morphism corresponds to a
nondegenerate pairing if
$ X $ is a rational homology manifold), because
$ \varphi_{f}\bQ_{X}|_{L_{x}} $ vanishes,
see also (1.3) below.

Let
$ b_{\lambda}^{j}(F_{x}) $ denote the rank of
$ H^{j}(F_{x},\bC)_{\lambda} \,(= \Ker\,(T_{s}-\lambda)) $
for
$ \lambda \in \bC $.
Using Theorem (0.1), we can explicitly calculate it for
$ j \le n - 2 $ in the case of a divisor with simple normal
crossings outside a point as follows
(see (4.3) for the proof).

\medskip\noindent
{\bf 0.2.~Theorem.}
{\it With the notation and the assumption of {\rm (0.1),}
assume
$ X \setminus \{x\} $ is smooth,
$ Y \setminus \{x\} $ is a divisor with normal
crossings on
$ X \setminus \{x\} $, and the local irreducible components
$ Y_{i}\,(i = 1,\dots, m) $ of
$ Y_{\red} $ at
$ x $ are principal divisors having at most isolated
singularities at
$ x $.
Let
$ a_{i} $ be the multiplicity of
$ Y $ at the generic point of
$ Y_{i} $, and
$ d = \hbox{\rm GCD}(a_{1},\dots, a_{m}) $.
Assume
$ j \le n - 2 + \delta_{\lambda,1} $, where
$ \delta_{\lambda,1} = 1 $ if
$ \lambda = 1 $, and
$ 0 $ otherwise.
Then
$ H^{j}(F_{x},\bQ) $ is a pure Hodge structure of type
$ (j,j) ; $ in particular, the monodromy is
semisimple.
Furthermore, if
$ \lambda^{d} \ne 1 $, we have
$ b_{\lambda}^{j}(F_{x}) = 0 $, and if
$ \lambda^{d} = 1 $, then
$$
\aligned
& b_{\lambda}^{j}(F_{x}) =
\hbox{$ \binom {m-1}j $}\,\,\,
\text{for}\,\, j < n - 2 + \delta_{\lambda,1},
\\
& b_{\lambda}^{j}(F_{x}) \le
\hbox{$ \binom {m-1}j $}\,\,\,
\text{for}\,\, j = n - 2 + \delta_{\lambda,1}.
\endaligned
$$
Here the equality holds also for
$ j = n - 2 + \delta_{\lambda,1} $, if
$ Y_{I} := \mcap_{i\in I}Y_{i} $ is a rational homology
manifold for any subset
$ I $ of
$ \{1,\dots, m\} $ with
$ |I| \le n-1 $, where
$ Y_{\emptyset} = X $.
}

\medskip
The case
$ a_{i} = 1 $ for any
$ i $ was studied in [9], see also (4.4) below.
In the case where an embedded resolution of
$ (X,Y) $ can be obtained by one blow-up with a point
center (e.g. an equisingular deformation of the affine
cone of a divisor with simple normal crossings on a
smooth projective variety),
we have a more precise statement as follows
(see (4.5) for the proof).

\medskip\noindent
{\bf 0.3.~Theorem.}
{\it With the notation and assumptions of {\rm (0.1)}, let
$ \pi : \tilde{X} \to X $ be the blow-up of
$ X $ with center
$ x $,
and assume that
$ \tilde{X} $ and the exceptional divisor
$ E := \pi^{-1}(x) $ are smooth and the total transform
$ \tilde{Y} := \pi^{-1}(Y) $ is a divisor with normal
crossings.
Let
$ Y' $ be the proper transform of
$ Y $, and put
$ U = E \setminus Y' $.
Let
$ e $ be the multiplicity of
$ \tilde{Y} $ along
$ U $.
Then the monodromy
$ T $ on
$ H^{j}(F_{x},\bQ) $ is semisimple for any
$ j $, and
$ H^{j}(F_{x},\bQ) $ is of type
$ (j,j) $ for
$ j < n $.
Furthermore, if
$ \lambda^{e} \ne 1 $, we have
$ b_{\lambda}^{j}(F_{x}) = 0 $ for any
$ j $, and if
$ \lambda^{e} = 1 $, then
$$
\aligned
& \chi_{\lambda}(F_{x}) \,(:=
\msum_{0\le j\le n}(-1)^{j}b_{\lambda}^{j}(F_{x})) =
\chi(U),
\\
& b_{\lambda}^{j}(F_{x}) =
\cases
\hbox{$ \binom {m-1}j $}\qquad &\text{if}\,\,
j < n ,\, \lambda^{d} = 1
\\
\,\,\,\,\,\, 0 &\text{if}\,\,
j < n ,\, \lambda^{d} \ne 1,
\endcases
\\
& b_{\lambda}^{n}(F_{x}) =
\cases
(-1)^{n}\chi(U)+
\hbox{$ \binom {m-2}{n-1} $} &\text{if}\,\,
\lambda^{d} = 1
\\
(-1)^{n}\chi(U) &\text{if}\,\,
\lambda^{d} \ne 1.
\endcases
\endaligned
$$
}

\medskip
This gives a generalization of formulas in [5], [15] for
a generic central arrangement with
$ a_{i} = 1 $, see (4.6) below.
If
$ X $ is smooth (i.e.
if
$ (X,x) = (\bC^{n+1},0)) $, then the assumption of (0.3)
is equivalent to that the union of the divisors defined by
the lowest degree part of a defining equation
$ f_{j} $ of
$ Y_{j} $ is a reduced divisor with normal crossings on
$ \bP^{n} $, and we have
$ e = \msum_{j} a_{j}d_{j} $ where
$ d_{j} $ is the degree of the lowest degree part of
$ f_{j} $; in particular,
$ d $ divides
$ e $.
We can calculate
$ \chi(U) $ explicitly in this case, see (4.6).

Let
$ T $ denote the monodromy of
$ \varphi_{f}\bQ_{X} $ with the Jordan decomposition
$ T = T_{u}T_{s} $.
For a complex number
$ \lambda $, set
$ \varphi_{f,\lambda}\bC_{X} = \Ker\,(T_{s}-\lambda)\subset
\varphi_{f}\bC_{X} $ (in the abelian category of shifted
perverse sheaves), and
$ N = \log T_{u} $.
As an application of Theorem (0.1), we show

\medskip\noindent
{\bf 0.4.~Theorem.}
{\it
With the notation and the assumption of {\rm (0.1)}, let
$ j $ be a positive integer
$ < n $.
Assume the monodromy of
$ \tH^{j}(F_{x},\bC)_{\lambda} $ has a Jordan block
of size
$ k $.
Then the action of
$ N^{k-1} $ on
$ \varphi_{f,\lambda}\bC_{X}|_{U\setminus \{x\}} $
is nonzero for any open neighborhood
$ U $ of
$ x $.
Furthermore, there exist points
$ y_{i} \,(\ne x) $ sufficiently near
$ x $ for
$ i \le j $ such that the monodromy of
$ \tH^{i}(F_{y_{i}},\bC)_{\lambda} $ has a Jordan block of
size
$ k_{i} $ and
$ \msum_{i\le j}k_{i} \ge k $, where we set
$ k_{i} = 0 $ if
$ \tH^{i}(F_{y},\bC)_{\lambda} = 0 $ for
$ y \ne x $.
}

\medskip
This is a refinement of Cor. 6.1.7 in [8].
There is an example such that the monodromy at degree
$ n - 1 $ is not semisimple at
$ x $, but is semisimple outside
$ x $, see Appendix.
Note that the support of the image of
$ N^{k} $ in
$ \psi_{f}\bQ $ (resp. in
$ \varphi_{f,1}\bQ $) as shifted perverse sheaves has
dimension
$ \le n - k $ (resp.
$ \le n - k - 1 $), see e.g. [10].
In the case
$ \dim \supp \varphi_{f,\lambda}\bC_{X} = r $, we have
$ \cH^{j}\varphi_{f,\lambda}\bC_{X} = 0 $ for
$ j < n - r $ (see (2.1.2) below), and the conclusion of
Theorem (0.4) for
$ j = n - r $ means that the monodromy of
$ \tH^{n-r}(F_{y},\bC)_{\lambda} $ has a Jordan block of
size
$ m $ for any point
$ y $ of a connected component of
$ L_{x} \cap \supp \varphi_{f,\lambda}\bC_{X} $
(considering the subsheaf of
$ \cH^{n-r}\varphi_{f,\lambda}\bC_{X} $ defined by the image
of
$ N^{k-1} $ and using (3.5) below).
In particular, we get

\medskip\noindent
{\bf 0.5.~Corollary.}
{\it If
$ \dim \supp \varphi_{f,\lambda}\bC_{X} = r $ {\rm (}e.g.
if
$ \dim \Sing f = r) $ and the monodromy of
$ \tH^{n-r}(F_{y},\bC)_{\lambda} $ for one point
$ y $ of each connected component of
$ L_{x} \cap \supp \varphi_{f,\lambda}\bC_{X} $ is semisimple,
then so is that of
$ \tH^{n-r}(F_{x},\bC)_{\lambda} $.
}

\medskip
For the lowest degree part we have a more precise description
of
$ \cH^{n-r}\varphi_{f,\lambda}\bC_{X} $, see (3.5) below.

In Sect.~1 we review the theory of nearby and vanishing cycles.
In Sect.~2 we calculate the cohomology of some sheaf complexes
on the link of a point.
In Sect.~3 we prove Theorems (0.1) and (0.4).
In Sect.~4 we treat the case of simple normal crossings outside
a point, and prove Theorems (0.2) and (0.3).
In Appendix we give a nontrivial example for Theorem (0.4).

\bigskip\bigskip
\centerline{{\bf 1. Vanishing Cycles}}

\bigskip\noindent
{\bf 1.1.~Nearby and vanishing cycles.}
Let
$ f $ be a nonconstant holomorphic function on a connected
complex analytic space
$ X $.
Assume
$ \bQ_{X}[n+1] $ is a perverse sheaf in the sense of [2]
(in particular,
$ \dim X = n + 1 $).
This is satisfied if
$ X $ is a locally complete intersection,
see e.g. [8], Th. 5.1.19.
(Indeed, if
$ X $ is defined locally by a regular sequence
$ g_{1},\dots, g_{r} $ on a smooth space
$ Z $, we can show the acyclicity (except for one degree)
of the algebraic local cohomology of
$ \cO_{Z} $ along
$ X $ by using the inductive limit of the Koszul complex of
$ g_{1}^{m},\dots, g_{r}^{m} $ for
$ m \to \infty $, see also (1.6) below.)

Let
$ A $ be a field of characteristic
$ 0 $ (e.g.
$ A = \bQ $ or
$ \bC $).
We denote by
$ \psi_{f}A_{X}, \varphi_{f}A_{X} $ the
nearby and vanishing cycle complexes on
$ Y := f^{-1}(0) $,
see [7].
It is well known that
$ \psi_{f}A_{X}[n] $ and
$ \varphi_{f}A_{X}[n] $ are perverse sheaves.
(This follows, for example, from [12], [13], see also [3].)
We have the action of the semisimple part
$ T_{s} $ of the monodromy
$ T $ on the shifted perverse sheaves.
For
$ \lambda \in A $, let
$$
\psi_{f,\lambda}A_{X} = \Ker\,(T_{s}-\lambda) \subset
\psi_{f}A_{X}\,\,\,
\text{(similarly for}\,\,\varphi_{f,\lambda}A_{X}).
$$
By definition of vanishing cycles, we have
$$
\psi_{f,\lambda}A_{X} = \varphi_{f,\lambda}A_{X}
\quad\text{for}\,\,\lambda \ne 1.
$$

If
$ A $ is algebraically closed, we have the decompositions
$$
\psi_{f}A_{X} = \mopls_{\lambda} \psi_{f,\lambda}A_{X},
\quad \varphi_{f}A_{X} =
\mopls_{\lambda} \varphi_{f,\lambda}A_{X},
$$
In general, we have
$$
\psi_{f}A_{X} = \psi_{f,1}A_{X} \oplus \psi_{f,\ne 1}A_{X},
\quad \varphi_{f}A_{X} = \varphi_{f,1}A_{X} \oplus
\varphi_{f,\ne 1}A_{X},
$$
where
$ \psi_{f,\ne 1} $,
$ \varphi_{f,\ne 1} $ denote the non unipotent monodromy part,
and
$ \psi_{f,\ne 1} = \varphi_{f,\ne 1} $.

For
$ x \in Y $,
we have isomorphisms
$$
H^{j}(F_{x},A)_{\lambda} = \cH^{j}(\psi_{f,\lambda}A_{X})_{x},
\quad \tH^{j}(F_{x},A)_{\lambda} =
\cH^{j}(\varphi_{f,\lambda}A_{X})_{x}.
\leqno(1.1.1)
$$
Here
$ F_{x} $ denotes the Milnor fiber as in the
introduction, and
$ H^{j}(F_{x},A)_{\lambda} $ is the
$ \lambda $-eigenspace as above.
By [16], [17], we have a canonical mixed Hodge structure
on these groups (which coincides with the one in [20] for
the isolated singularity case), see also [14].

\medskip\noindent
{\bf 1.2.~Cohomology with compact supports.}
It is known that there is a proper continuous map
$ \rho : X_{c} \to Y $ such that
$ \psi_{f}A = \boR\rho_{*}A $,
where
$ X_{c} = f^{-1}(c) $ for
$ c \ne 0 $ sufficiently small.
This can be constructed by using a resolution of singularities.
Let
$ i : \{x\} \to Y $ denote the
inclusion morphism.
Then for a sufficiently small open ball
$ B_{x} $ around
$ x $,
we have a commutative diagram
$$
\CD
H_{c}^{k}(F_{x},A) @= H_{c}^{k}(B_{x}\cap Y,\psi_{f}A) @=
H^{k}i^{!}\psi_{f}A
\\
@VV{\beta_{F}}V @VV{\beta_{B}}V @VV{\beta_{\psi}}V
\\
H^{k}(F_{x},A) @= H^{k}(B_{x}\cap Y,\psi_{f}A) @=
H^{k}i^{*}\psi_{f}A,
\endCD
\leqno(1.2.1)
$$
where the horizontal morphisms are canonical isomorphisms,
the first two vertical morphisms
$ \beta_{F} $,
$ \beta_{B} $ are natural morphisms, and
$ \beta_{\psi} $ is induced by the natural morphism
$ i^{!} \to i^{*} $.
By (1.2.1),
$ \beta_{F} $ will be identified with
$ \beta_{\psi} $.

\medskip\noindent
{\bf 1.3.~Unipotent monodromy part.}
We have morphisms of perverse sheaves (compatible with mixed
Hodge modules [16])
$$
\can : \psi_{f,1}A \to \varphi_{f,1}A,\quad
\Var : \varphi_{f,1}A(1) \to \psi_{f,1}A,
$$
whose compositions coincide with
$ N $ on
$ \psi_{f,1}A, \varphi_{f,1}A $.
If
$ n \ge 1 $ and
$ X $ is a rational homology manifold at
$ x $, then they induce isomorphisms
$$
\aligned
\can
&: H^{n}i^{*}\psi_{f,1}A \simto H^{n}i^{*}\varphi_{f,1}A,
\\
\Var
&: H^{n}i^{!}\varphi_{f,1}A(1) \simto H^{n}i^{!}\psi_{f,1}A,
\endaligned
\leqno(1.3.1)
$$
because the mapping cone of
$ \Var $ is
$ \boR\Gamma_{Y}A_{X}(1)[2] $ and
$ \boR\Gamma_{\{x\}}\boR\Gamma_{Y}A_{X} =
\boR\Gamma_{\{x\}}A_{X} $.

By the isomorphisms of (1.3.1), the morphism
$$
\beta_{F,1} : H_{c}^{n}(F_{x},A)_{1} \to H^{n}(F_{x},A)_{1},
\leqno(1.3.2)
$$
which is the restriction of
$ \beta_{F} $, can be identified with the composition of
$ N $ and
$$
\beta_{\varphi,1} : H^{n}i^{!}\varphi_{f,1}A \to
H^{n}i^{*}\varphi_{f,1}A,
\leqno(1.3.3)
$$
which is induced by the natural morphism
$ i^{!} \to i^{*} $.
Indeed, using
$ \can\scirc\Var = N $ together with the commutativity of the
natural morphism
$ i^{!} \to i^{*} $ with
$ \can $,
$ \Var $, we get a commutative diagram
$$
\CD
H^{n}i^{!}\psi_{f,1}A @>{\beta_{\psi,1}}>>
H^{n}i^{*}\psi_{f,1}A
\\
@AA{\Var}A @VV{\can}V
\\
H^{n}i^{!}\varphi_{f,1}A(1) @>{N\beta_{\varphi,1}}>>
H^{n}i^{*}\varphi_{f,1}A
\endCD
\leqno(1.3.4)
$$
where the vertical morphisms are isomorphisms.
Note that the morphism
$ \beta_{\varphi,1} $ in (1.3.3) is an isomorphism in the
isolated singularity case, because
$ \supp \varphi_{f}A = \{x\} $.

In Theorem (0.1),
$ \beta_{\varphi,1} $ in (1.3.3) is identified with a morphism
$ H^{n}_{c}(F_{x},A)_{1}(-1) \to H^{n}(F_{x},A)_{1} $ by using
the isomorphisms of (1.2.1) and (1.3.1).
For the non unipotent monodromy part, we have
$ \beta_{\psi,\ne 1} = \beta_{\varphi,\ne 1} $,
because
$ \psi_{f,\ne 1} = \varphi_{f,\ne 1} $.

\medskip\noindent
{\bf 1.4.~Normal crossing case.}
Assume that
$ Y := f^{-1}(0) $ is a divisor with normal crossings on a
complex manifold
$ X $ whose irreducible components
$ Y_{1},\dots, Y_{m} $ are smooth.
Let
$$
\cF_{\lambda} = \psi_{f,\lambda}\bC_{X}[n].
$$
Since
$ \psi_{f,\lambda}\bC_{X} \oplus \psi_{f,\olam}\bC_{X} $
underlies a mixed Hodge Module,
$ \cF_{\lambda} $ has the weight filtration
$ W $ which is the monodromy filtration shifted by
$ n = \dim Y $,
i.e.
$$
N^{k} : \Gr_{n+k}^{W}\cF_{\lambda} \simto
\Gr_{n-k}^{W}\cF_{\lambda}.
\leqno(1.4.1)
$$
Let
$ P\Gr_{n+k}^{W}\cF_{\lambda} $ denote the
$ N $-primitive part, which is defined by
$ \Ker\, N^{k+1} \subset \Gr_{n+k}^{W}\cF_{\lambda} $ for
$ k \ge 0 $, and is zero otherwise.
By (1.4.1) we have the primitive decomposition
$$
\Gr_{j}^{W}\cF_{\lambda} =
\mopls_{k\ge 0}N^{k}P\Gr_{j+2k}^{W}\cF_{\lambda}(k).
\leqno(1.4.2)
$$

Let
$ a_{j} $ be the multiplicity of
$ f $ along
$ Y_{j} $, and put
$$
J(\lambda) = \{j : \lambda^{a_{j}} = 1\}.
$$
Let
$ d = \hbox{\rm GCD}(a_{1},\dots,a_{m}) $.
Then
$$
J(\lambda) = \{1,\dots, m\}\,\,\,
\text{if and only if
$ \lambda^{d} = 1 $.}
\leqno(1.4.3)
$$
For
$ I \subset J(\lambda) $,
let
$$
Y_{I} = \mcap_{j\in I} Y_{j},\quad
U_{I} = Y_{I} \setminus \mcup_{j\notin J(\lambda)} Y_{j},
$$
with the inclusion morphism
$ j_{I} : U_{I} \to Y_{I} $.
By [17], 3.3, we see that the primitive part
$ P\Gr_{n+k}^{W}\cF_{\lambda} $ is the direct sum of
$$
(j_{I})_{!}\cF_{\lambda,I}(-k)[n-k] =
\bold{R}(j_{I})_{*}\cF_{\lambda,I}(-k)[n-k]
\leqno(1.4.4)
$$
over
$ I \subset J(\lambda) $ with
$ |I| = k + 1 $,
where
$ \cF_{\lambda,I} $ is a local system of rank
$ 1 $ on
$ U_{I} $.
Furthermore, the monodromy of
$ \cF_{\lambda,I} $ around
$ Y_{j} \,(j \notin J(\lambda)) $ is given by the
multiplication by
$ \lambda^{-a_{j}} $ so that (1.4.4) holds.

If each
$ Y_{j} $ is a principal divisor defined by a reduced
equation
$ f_{j} $ and
$ f = \prod_{j} f_{j}^{a_{j}} $, then the
$ \cF_{\lambda,I} $ are the restrictions of
$ \cF_{\lambda,\emptyset} $ on
$ U_{\emptyset} $ which is defined by
$ \motim_{j}f_{j}^{*}L_{j}$ where
$ L_{j} $ is a local system on
$ \bC^{*} $ with monodromy
$ \lambda^{-a_{j}} $ for
$ j \notin J(\lambda) $.
This can be verified by reducing to the case where the
$ a_{i} $ are independent of
$ i $, and using the compatibility of the nearby cycle
functor with the direct image under a proper morphism.
Indeed, setting
$ c_{j} = \hbox{\rm LCM}(a_{1},\dots,a_{m})/a_{j} $,
we have a ramified covering of
$ X $ defined by
$$
\{(x,t_{1},\dots,t_{m})\in X\times\bC^{m}:
f_{j}(x)=t_{j}^{c_{j}}\,\,\text{for any}\,\,j\}.
\leqno(1.4.5)
$$

For the vanishing cycle
$ \varphi_{f,1}\bC_{X}[n] $ with
$ \lambda = 1 $,
the weight filtration is the monodromy filtration shifted by
$ n + 1 $.
For the
$ N $-primitive part
$ P\Gr_{n+1+k}^{W}\varphi_{f,1}\bC_{X}[n] $,
we have
$$
P\Gr_{n+1+k}^{W}\psi_{f,1}\bC_{X}[n] =
P\Gr_{n+1+k}^{W}\varphi_{f,1}\bC_{X}[n]\quad\text{for}\,\,k \ge 0,
$$
because
$ \varphi_{f,1}\bC_{X}[n] $ can be identified with
$ \Im\, N \subset \psi_{f,1}\bC_{X}[n] $.

\medskip\noindent
{\bf 1.5.~Weight spectral sequence.}
$ $ Let
$ \pi : (X',Y') \to (X,Y) $ be an embedded resolution such that
$ Y' := \pi^{-1}(Y) $ and
$ E := \pi^{-1}(x) $ are divisors with simple normal crossings.
Let
$ E' $ be the closure of
$ Y' \setminus E $, and put
$ U = E \setminus E' $ with the inclusion
$ j' : U \to E $.
Let
$ f' = f\pi $.
Then by [4], 4.2, the canonical morphism
$$
\psi_{f',\lambda}\bC_{X}|_{E'} \to \bold{R}j'_{*}
(\psi_{f',\lambda}\bC_{X}|_{U})
\leqno(1.5.1)
$$
is a quasi-isomorphism.
(This easily follows from [17], 3.3.)
Since the nearby cycle functor commutes with the direct
image under a proper morphism, we get canonical isomorphisms
(compatible with
$ T) $
$$
H^{i}(F_{x},\bC)_{\lambda} =
H^{i}(E,\psi_{f,\lambda}\bC_{X}|_{E}) =
H^{i}(U,\psi_{f,\lambda}\bC_{X}|_{U}).
\leqno(1.5.2)
$$

Let
$ Y_{1},\dots, Y_{m} $ denote the irreducible components of
$ Y' $ (which are assumed to be smooth).
We may assume that
$ Y_{1},\dots, Y_{r} $ are the irreducible components of
$ E = \pi^{-1}(x) $.
Let
$ Y_{I}, U_{I}, \cF_{\lambda,I} $ be as in (1.4).
For
$ I \subset \{1,\dots, m\} $, let
$ s(I) = |I \cap \{1,\dots, r\}| - 1 $.
Then we have the weight spectral sequence
$$
E_{1}^{-k,j+k} = \mopls_{I,a}
H^{j-|I|+1}(U_{I},\cF_{\lambda,I}(a+1-|I|)) \Rightarrow
H^{j}(F_{x},\bC)_{\lambda},
\leqno(1.5.3)
$$
where the summation is taken over
$ I \,(\ne \emptyset) \subset J(\lambda ), 0 \le a \le s(I) $
such that
$ |I| - 1 -2a = k $.
Indeed,
$ \boR j'_{*} $ is a
$ t $-exact functor [2], and
$ (j_{I})_{!}\cF_{\lambda,I}(a+1-|I|)[n+1-|I|] $
comes from the graded pieces of the weight filtration on
$$
\bold{R}j'_{*}((j_{I'})_{!}\cF_{\lambda,I'}(a+1-|I'|)
[n+1-|I'|]|_{U})
$$
for
$ I' := I \cap \{1,\dots, r\} $.
Here we may assume essentially that
$ \cF_{\lambda,I'} $ is a constant sheaf
(where the assertion is well-known [6])
because it is of normal crossing type, see [17], 3.1.
The range of
$ a $ comes from the symmetry of the weight filtration
(1.4.1) which is related to
$ \cF_{\lambda,I'} $ because we consider it on
$ U $.

The spectral sequence (1.5.3) degenerates at
$ E_{2} $,
because
$ E_{1}^{-k,j+k} $ is pure of weight
$ j + k $.

\medskip\noindent
{\bf 1.6.~Remark.}
If
$ \bQ_{X}[n+1] $ is a perverse sheaf, then
$ \bQ_{Y}[n] $ is a perverse sheaf for any locally principal
divisor
$ Y $ on
$ X $.
Indeed, we have locally a distinguished triangle
$$
\bQ_{Y}[n] \to \psi_{f}\bQ_{X}[n] \to \varphi_{f}
\bQ_{X}[n] \buildrel{+1}\over\to,
\leqno(1.6.1)
$$
by the definition of
$ \varphi_{f} $,
where
$ f $ is a local equation of
$ Y $.
This implies
$ {}^{p}\cH^{j}(\bQ_{Y}[n]) = 0 $ except for
$ j = 0, 1 $,
where
$ {}^{p}\cH^{j} $ denotes the perverse cohomology functor [2].
Furthermore, the vanishing of
$ {}^{p}\cH^{j}(\bQ_{Y}[n]) $ for
$ j > 0 $ is clear by the definition of semi-perversity.
(In general, a sheaf complex
$ \cF $ is called semi-perverse if
$ \dim \supp \cH^{-i}\cF \le i $ for any
$ i $, see loc.~cit.)

\medskip\noindent
{\bf 1.7.~Wang sequence.}
Let
$ f $ be a holomorphic function on an analytic space
$ X $.
Let
$ L_{X,x} $ be the link of
$ x $ in
$ X $.
Then we have the Wang sequence
$$
H^{j}(L_{X,x}\setminus Y,\bQ) \to
H^{j}(F_{x},\bQ)_{1} \buildrel{N}\over\to
H^{j}(F_{x},\bQ)_{1}(-1) \to
H^{j+1}(L_{X,x}\setminus Y,\bQ).
$$
In the category of mixed Hodge structures, this follows from
$$
i^{\prime *}j'_{*}\bQ_{X} =
C(j'_{!}\bQ_{X} \to j'_{*}\bQ_{X}) =
C(N : \psi_{f,1}\bQ_{X} \to \psi_{f,1}\bQ_{X}(-1))[-1],
$$
where
$ i' : Y \to X, j' : X \setminus Y \to X $ are the
inclusion morphisms,
see e.g. [17], 2.23 for the second isomorphism.
(Here
$ \bQ_{X} $ can be defined locally in the derived category
of mixed Hodge Modules, using an embedding into a smooth
space.)

\bigskip\bigskip
\centerline{{\bf 2. Cohomology of Link with Coefficients}}

\bigskip\noindent
{\bf 2.1.~Localization sequence.}
Let
$ \cF $ be a perverse sheaf on
$ Y $ in the sense of [2].
In particular,

\bigskip\noindent
(2.1.1)\qquad\qquad
$ \dim \supp \cH^{-k}\cF \le k $,

\bigskip\noindent
(2.1.2)\qquad\qquad
$ \cH^{-r}\cF = 0 $ for
$ r > \dim \supp \cF $.

\bigskip
Let
$ i : \{x\} \to Y $ and
$ j : U := Y \setminus \{x\} \to Y $ denote the
inclusions.
Let
$ L_{x} $ be the intersection of a sufficiently small sphere
around
$ x $ with
$ Y $.
Then
$$
\boH^{k}(L_{x},\cF|_{L_{x}}) = H^{k}i^{*}j_{*}j^{*}\cF,
\leqno(2.1.3)
$$
and we get a long exact sequence
$$
\to H_{\{x\}}^{k}\cF \to H^{k}\cF_{x} \to
\boH^{k}(L_{x},\cF|_{L_{x}}) \to H_{\{x\}}^{k+1}\cF \to
$$
induced by the distinguished triangle
$$
\boR\Gamma_{\{x\}}\cF \to \cF_{x} \to
\boR\Gamma (L_{x},\cF|_{L_{x}}) \buildrel{+1}\over\to
$$
which is identified with
$ i^{!}\cF \to i^{*}\cF \to i^{*}j_{*}j^{*}\cF
\buildrel{+1}\over\to $ (because
$ i_{*}i^{!} = \boR\Gamma_{\{x\}} $).

Let
$ \bD $ denote the functor assigning the dual.
Since
$ \bD i^{!} = i^{*}\bD $,
and
$ \bD\cF $ is a perverse sheaf, we get
$$
H_{\{x\}}^{k}\cF = 0\quad\text{for}\,\,k < 0.
\leqno(2.1.4)
$$
Indeed, (2.1.4) is equivalent to the (dual)
semi perversity of
$ \cF $ (see [2]) assuming the perversity of the
restriction of
$ \cF $ to the complement of
$ x $.

\medskip\noindent
{\bf 2.2.~Leray spectral sequence.}
Let
$ \cF $ be a complex of sheaves with constructible cohomology on
$ Y $.
There is a Leray-type spectral sequence
$$
E_{2}^{p,q} = H^{p}(L_{x},\cH^{q}\cF|_{L_{x}}) \Rightarrow
\boH^{p+q}(L_{x},\cF|_{L_{x}})
\leqno(2.2.1)
$$
induced by the filtration
$ \tau $ on
$ \cF $, see [6].
By (2.1.3) this is compatible with mixed Hodge structure
(using a
$ t $-structure in [17], 4.6) if
$ \cF $ underlies a complex of mixed Hodge modules.
The calculation of (2.2.1) is not necessarily easy.
One problem is that
$ \cH^{q}\cF $ is a constructible sheaf and not a local system,
and some times we have to use the spectral sequence associated
to a stratification, which is a special case of (2.3.1) below,
to calculate its cohomology.
Actually this spectral sequence can be formulated for a complex
as below, and we do not have to use spectral sequences twice
if we can calculate the
$ E_{1} $-term of (2.3.1).
But the calculation of
$ d_{r} $ is still nontrivial.

\medskip\noindent
{\bf 2.3.~Spectral sequence associated to a stratification.}
Let
$ \cF $ be as above, and let
$ \{Y_{k}\} $ be a stratification of
$ Y $ compatible with
$ \cF $,
where the
$ Y_{k} $ are locally closed analytic subspaces of
$ Y $ with pure dimension
$ k $ such that the restriction of
$ \cH^{j}\cF $ to
$ Y_{k} $ is a local system, and
$ \oY_{k} \setminus Y_{k} $ is the disjoint union of
$ Y_{i} \,(i < k) $.
Put
$ U_{k} = Y \setminus \oY_{k-1} $.
Then, for each $k$, there is a subcomplex of
$ \cF $ whose restriction to
$ U_{k} $ coincides with
$ \cF|_{U_{k}} $ and whose restriction to
$ \oY_{k-1} $ vanishes (i.e. it is the direct image with
proper supports by
$ U_{k} \to Y $).
Such complexes form a decreasing filtration of
$ \cF $ whose graded pieces are (the direct images with
proper supports of) the restrictions of
$ \cF $ to the
$ Y_{k} $.
So they induce the spectral sequence associated to the
stratification
$$
E_{1}^{p,q} =
\boH_{c}^{p+q}(L_{x}\cap Y_{p}, \cF|_{L_{x}\cap Y_{p}})
\Rightarrow
\boH^{p+q}(L_{x},\cF|_{L_{x}}).
\leqno(2.3.1)
$$
By (2.1.3) this is also compatible with mixed Hodge structure
(using the quasi-filtration in [16], 5.2.17).

\medskip\noindent
{\bf 2.4.~Weight spectral sequence.}
Let
$ \cF $ be a perverse sheaf underlying a mixed Hodge Module,
and
$ W $ be the weight filtration.
Then, as in [6],
$ W $ induces a spectral sequence
$$
E_{1}^{-k,j+k} = \boH^{j}(L_{x},\Gr_{k}^{W}
\cF|_{L_{x}}) \Rightarrow \boH^{j}(L_{x},\cF|_{L_{x}}),
\leqno(2.4.1)
$$
which is called the (generalized) weight spectral sequence.
(We can use Verdier's theory of spectral objects,
see [2] and also [16], 5.2.18.)
By (2.1.3) this is compatible with mixed Hodge structure, but
does not necessarily degenerate at
$ E_{2}, $ because
$ E_{1}^{-k,j+k} $ is not pure of weight
$ j + k $ in general.
It is not easy to calculate this spectral sequence explicitly
except for some special cases, see e.g. (4.2) below.

If
$ X \setminus \{x\} $ is smooth and
$ Y \setminus \{x\} $ is a divisor with simple normal
crossings, then the
$ E_{1} $-complex has a structure of double complex whose
differentials are induced by the Cech restriction morphism
and the co-Cech Gysin morphism, see e.g. [20].
Indeed, the differential
$ d_{1} $ is induced by the extension class between the
graded pieces of the perverse sheaves, and the assertion
can be verified by using locally a ramified covering as
in (1.4.5) and reducing to the case where the irreducible
components of
$ Y \setminus \{x\} $ have the constant multiplicity.

\bigskip\bigskip
\centerline{{\bf 3. Proofs of Theorems (0.1) and (0.4)}}

\bigskip\noindent
{\bf 3.1.~Proof of Theorem (0.1).} Applying (2.1) to
$ \cF = \varphi_{f}\bQ_{X}[n] $,
the assertion follows from (1.1--3) and (2.1).

\medskip\noindent
{\bf 3.2.~Proof of Theorem (0.4).}
The first assertion follows from (2.1) applied to
$ \Im\,N^{k-1} \subset \varphi_{f}\bC_{X} $ (defined in the
abelian category of shifted perverse sheaves).
Indeed, factorizing
$ N^{k-1} : \varphi_{f}\bC_{X} \to \varphi_{f}
\bC_{X}(1-k) $ by
$ \Im\,N^{k-1} $,
we see that
$ \Im\,N^{k-1} \ne 0 $ on a neighborhood of
$ x $.
The remaining assertion is clear by (2.2.1).
Indeed, if any Jordan block of the monodromy on
$ \tH^{i}(F_{y},\bC)_{\lambda} $ has size at most
$ k_{i} $, then
$ N^{k_{i}} = 0 $ on
$ H^{j-i}(L_{x},\cH^{i}\varphi_{f,\lambda}\bC_{X}|_{L_{x}}) $,
and
$ N^{k}(\tH^{j}(F_{x},\bC)_{\lambda}) = 0 $ for
$ k = \msum_{i\le j}k_{i} $ by Theorem (0.1) together with
(2.2.1), because
$ N^{k_{i}}(\Gr_{G}^{j-i}\tH^{j}(F_{x},\bC)_{\lambda}) = 0 $
where
$ G $ is the filtration associated to the spectral sequence
(2.2.1).

\medskip\noindent
{\bf 3.3.~One-dimensional singular locus case.} If
$ \Sigma_{\lambda}:=\supp\varphi_{f,\lambda}\bC_{X} $
is
$ 1 $-dimensional (e.g. if
$\Sing f $ is
$ 1 $-dimensional), let
$ \Sigma_{\lambda,i} $ be the local irreducible components of
$ \Sigma_{\lambda} $ at
$ x $,
and take
$ x_{i} \in \Sigma_{\lambda,i} \cap L_{x} $.
Then
$ H^{j}\varphi_{f,\lambda}\bC_{X} = 0 $ for
$ j < n - 1 $,
and
$$
\boH^{n-1}(L_{x},\varphi_{f,\lambda}\bC_{X}|_{L_{x}})
=\mopls_{i}(\tH^{n-1}(F_{x_{i}},\bC)_{\lambda})^{\tau_{i}},
\leqno(3.3.1)
$$
where
$ \tau_{i} $ denotes the monodromy of the local system
on
$ \Sigma_{\lambda,i} \cap L_{x} $ (which is called the
vertical monodromy in [18], [19]).
However, for a given element of
$ \mopls_{i}(\tH^{n-1}(F_{x_{i}},\bC)_{\lambda})^{\tau_{i}} $,
it is not easy to determine whether it comes from
$ \tH^{n-1}(F_{x},\bC)_{\lambda} $ or not.
Note that
$ K_{x} = \Ker\,\beta_{\varphi} $ does not vanish in general.
For example, if
$ X $ is smooth and
$ Y $ is a reduced divisor with normal crossings it is
well-known (see e.g. [20]) that the Milnor fiber is
homotopy equivalent to a real torus of dimension
$ m - 1 $ where
$ m $ is the multiplicity of
$ Y $ at the point.
In the case
$ f = xyz $ and
$ n = 2 $, we have
$ \dim H^{1}(F_{x},\bC)_{1} = 2 $ and
$ \dim \boH^{1}(L_{x},\varphi_{f,1}\bC_{X}|_{L_{x}}) = 3 $,
see also [18], [19].

\medskip\noindent
{\bf 3.4.~Remark.}
There are examples such that the monodromy of
$ \tH^{n-1}(F_{x},\bQ) $ is semisimple, but that of
$ \tH^{n-1}(F_{y},\bQ) $ for
$ y $ sufficiently near
$ x $ has a Jordan block of size
$ n $ (this implies that the converse of Theorem (0.4)
does not hold).
For example, consider a germ of
$ (n-1) $-dimensional hypersurface
$ (Y,x) $ with isolated singularity whose Milnor monodromy has
a Jordan block of size
$ n $, take a projective compactification
$ Z $ of
$ Y $ in
$ \bP^{n} $ such that
$ Z \setminus \{x\} $ is smooth (using finite determinacy of
isolated singularity), and then take
$ f : \bC^{n+1} \to \bC $ to be a defining equation of
$ Z $.

\medskip\noindent
{\bf 3.5.~Lowest degree term.}
Assume
$ \Sigma_{\lambda} := \supp \varphi_{f,\lambda}\bC_{X} $ is
$ r $-dimensional (e.g.
$ \Sing f $ is
$ r $-dimensional).
Let
$ \Sigma_{\lambda}^{1} $ be an
$ (r-1) $-dimensional Zariski-locally closed smooth analytic
subspace of
$ \Sigma_{\lambda} $ such that
$ \Sigma_{\lambda}^{0} := \Sigma_{\lambda} \setminus
\oSig{}_{\lambda}^{1} $ is smooth (where
$ \oSig{}_{\lambda}^{1} $ is the closure of
$ \Sigma_{\lambda}^{1} $) and the restrictions of
$ \cH^{j}\varphi_{f,\lambda}\bC_{X} $ to
$ \Sigma_{\lambda}^{0} $,
$ \Sigma_{\lambda}^{1} $ are local systems for any
$ j $.
Let
$ \oSig{}_{\lambda}^{2} = \oSig{}_{\lambda}^{1} \setminus
\Sigma_{\lambda}^{1} $,
$ U_{\lambda} = \Sigma_{\lambda}\setminus\oSig{}_{\lambda}^{2} $
with the inclusions
$ j' : \Sigma_{\lambda}^{0} \to U_{\lambda} $,
$ j'' : U_{\lambda} \to \Sigma_{\lambda} $.
Then

\bigskip\noindent
(3.5.1)\qquad\qquad
$ \cH^{n-r}\varphi_{f,\lambda}\bC_{X}|_{U_{\lambda}}
\subset j'_{*}j^{\prime *}(\cH^{n-r}\varphi_{f,\lambda}
\bC_{X}|_{U_{\lambda}}) $,

\bigskip\noindent
(3.5.2)\qquad\qquad
$ \cH^{n-r}\varphi_{f,\lambda}\bC_{X} = j''_{*}
(\cH^{n-r}\varphi_{f,\lambda}\bC_{X}|_{U_{\lambda}}) $.

\bigskip
Indeed, restricting to a subspace transversal to
$ \Sigma_{\lambda}^{1} $, (3.5.1) follows from the
$ 1 $-dimensional singular locus case, and furthermore,
the cokernel of the inclusion in (3.5.1) is given by
$ K_{x} $ in Theorem (0.1), see (3.3).
Similarly (3.5.2) follows from Theorem (0.1) by induction on
strata.

\bigskip\bigskip
\centerline{{\bf 4. Case of Simple Normal Crossings outside
a Point}}

\bigskip\noindent
{\bf 4.1.}
With the notation of (1.1), assume that
$ X \setminus \{x\} $ is smooth, and
$ Y \setminus \{x\} $ is a divisor with {\it simple}
normal crossings on
$ X \setminus \{x\} $.
Here simple means that each irreducible component of
$ Y \setminus \{x\} $ is smooth.
Assume further that the local irreducible components of
$ Y $ at
$ x $ are principal divisors.
Then, replacing $ X $ with a sufficiently small open
neighborhood of
$ x $ if necessary, there exist holomorphic functions
$ f_{i} : X \to \bC $ and positive integers
$ a_{i} $ for
$ i = 1,\dots, m $ such that
$ f = {f}_{1}^{{a}_{1}} \cdots {f}_{m}^{{a}_{m}} $ and each
$ Y_{i} := {f}_{i}^{-1}(0) $ has at most isolated
singularity at
$ x $, see also [9].
Here we assume
$ n \ge 2 $.
Let
$$
\cF_{\lambda} = \psi_{f,\lambda}\bC_{X}[n]|_{Y\setminus \{x\}}.
\leqno(4.1.1)
$$
Since
$ \psi_{f,\lambda}\bC_{X}[n] \oplus \psi_{f,\olam}\bC_{X}[n] $
underlies a mixed Hodge Module, we have the weight spectral
sequence (2.4.1).
Here
$ W $ is the monodromy filtration shifted by
$ n = \dim Y $,
and the
$ N $-primitive part
$ P\Gr_{n+k}^{W}\cF_{\lambda}|_{Y\setminus \{x\}} $ is
calculated as in (1.4).

We assume that
$ \bQ_{X}[n+1] $ is a perverse sheaf.
Since the intersection complex of
$ X $ is given by
$ \tau_{<0}\boR j'_{*}\bQ_{X\setminus \{x\}}[n+1] $
where
$ j' : X \setminus \{x\} \to X $ denotes the inclusion,
this condition is equivalent to
$$
\tH^{j}(L_{X,x},\bQ) = 0\quad \text{for}\,\,j < n,
\leqno(4.1.2)
$$
where
$ L_{X,x} $ is the link of
$ \{x\} $ in
$ X $.
This follows from the long exact sequence of perverse sheaves
associated to the distinguished triangle
$$
\bQ_{X}[n+1] \to
\tau_{<0}\boR j'_{*}\bQ_{X\setminus \{x\}}[n+1] \to
(\tau_{<0}\boR j'_{*}\bQ_{X\setminus \{x\}}/\bQ_{X})[n+1]
\buildrel{+1}\over\to,
\leqno(4.1.3)
$$
because
$ \cH^{j}(\boR j'_{*}\bQ_{X\setminus \{x\}})_{x} =
H^{j}(L_{X,x},\bQ) $.

\medskip\noindent
{\bf 4.2.~Proposition.}
{\it With the above notation and assumptions, let
$ \cF_{\lambda,I} $,
$ j_{I} $ and
$ d $ be as in {\rm (1.4)} with
$ Y_{j} $ replaced by
$ Y_{j} \setminus \{x\} $.
Then
$$
\boH^{i}(L_{x},(j_{I})_{!}\cF_{\lambda,I}[n-k]|_{L_{x}}) = 0
\quad\text{for}\,\,k - n < i < -1,
\leqno(4.2.1)
$$
where
$ k = |I| - 1 $.
For
$ i = k - n < -1 $, we have
$$
\boH^{k-n}(L_{x},(j_{I})_{!}\cF_{\lambda,I}[n-k]|_{L_{x}}) =
\cases
\bC &\text{if
$ \lambda^{d} = 1 $}
\\
0 &\text{if
$ \lambda^{d} \ne 1 $}.
\endcases
\leqno(4.2.2)
$$
}

\medskip\noindent
{\it Proof.}
We prove the assertion by induction on
$ |I| $.
If
$ I = \emptyset $,
we have
$ \cF_{\lambda,\emptyset} $ on
$ U_{\emptyset} $ as in (1.4).
We may assume
$ U_{\emptyset} \ne X $, because the assertion is clear by
(4.1.2) if
$ U_{\emptyset} = X $.
Let
$ B_{x} $ be a sufficiently small open ball around
$ x $.
By the cone theorem,
$ B_{x}\cap Y $ is homeomorphic to the topological cone of
$ \partial B_{x}\cap Y $ in a compatible way with a
given Whitney stratification of
$ Y $.
(This is proved by using a continuous vector field compatible
with the stratification as well-known.)
So we have
$$
\boH^{i}(L_{x},\boR (j_{\emptyset})_{*}
\cF_{\olam,\emptyset}|_{L_{x}}) =
\boH^{i}(B_{x}\cap U_{\emptyset},\cF_{\olam,\emptyset}).
$$
By duality, (4.2.1) is equivalent to the vanishing of these
groups for
$ n+1 < i < 2n+1 $.
(Note that the dual of the
$ \lambda $-eigenspace is the
$ \lambda^{-1} $-eigenspace.)
So we get the assertion in this case, using the
corresponding de Rham complex and the vanishing of the
higher cohomology of coherent sheaves on a smooth Stein space
$ B_{x}\cap U_{\emptyset} $
of dimension
$ n + 1 $.

If
$ I \ne \emptyset $,
take
$ j \in I $,
and let
$ I' = I \setminus \{j\} $.
By the exact sequence
$$
H^{i-1}(L_{x},(j_{I'})_{!}\cF_{\lambda,I'}) \to
H^{i-1}(L_{x},(j_{I})_{!}\cF_{\lambda,I}) \to
H_{c}^{i}(L_{x}\setminus Y_{j},(j_{I'})_{!}\cF_{\lambda,I'}),
$$
it is enough to show
$$
H_{c}^{i}(L_{x}\setminus Y_{j},(j_{I'})_{!}
\cF_{\lambda,I'}) = 0\quad \text{for}\,\,i < n - k.
$$
This is isomorphic to the dual of
$ H^{2n-2k+1-i}(L_{x}\setminus Y_{j},\bold{R}(j_{I'})_{*}
\cF_{\olam,I'}) $, because
$ \dim U_{I'} = n - k + 1 $.
So it is enough to show
$$
H^{i}(L_{x} \cap (Y_{I'} \setminus Y_{j}),\bold{R}
(j_{I'})_{*}\cF_{\olam,I'}) = 0\quad
\text{for}\,\,i > n - k + 1.
$$
For this, we may replace
$ L_{x} \cap (Y_{I'} \setminus Y_{j}) $ by
$ B_{x} \cap (U_{I'} \setminus Y_{j}) $ (using the cone
theorem).
Then we get the assertion by using the same argument as above,
because
$ B_{x} \cap (U_{I'} \setminus Y_{j}) $ is a smooth Stein
space of dimension
$ n - k + 1 $.

\medskip\noindent
{\bf 4.3.~Proof of Theorem (0.2).}
If
$ \lambda^{d} \ne 1 $, the assertion follows from (4.2).
So we may assume
$ \lambda^{d} = 1 $, i.e.
$ J(\lambda) = \{1,\dots, m\} $, see (1.4.3).
We define
$ K_{\lambda} $ to be a complex whose
$ j $-th component is
$$
\mopls_{|I|=j}H^{0}(L_{X,x}\cap Y_{I},\bC),
$$
where
$ Y_{\emptyset} = X $, and the differential is given by
the Cech restriction morphism.
Let
$ \sigma $ be the filtration as in [6], II, 1.4.7, and define
$$
\tK_{\lambda} =
\mopls_{i\ge 1} (\sigma_{\ge i}K_{\lambda})(1-i)[n+1].
\leqno(4.3.1)
$$

Let
$ E_{1,\lambda} $ denote the
$ E_{1} $-complex of the weight spectral sequence (2.4.1)
applied to (4.1.1).
For
$ -n \le j \le -2 $, we see that
$$
E_{1}^{-k,j+k} = K_{\lambda}^{j+n+1}
(\hbox{$ \frac {-j-k}2 $}),
\leqno(4.3.2)
$$
if
$ j + k $ is even and
$ |k| + j + n > 0 $, and it is zero otherwise,
using (4.2) and the primitive decomposition (1.4.2).
So we get
$$
\sigma_{\le -2}\tK_{\lambda} = \sigma_{\le -2}E_{1,\lambda},
\leqno(4.3.3)
$$
and
$ \sigma_{\le -1}\tK_{\lambda} $ is a quotient complex of
$ \sigma_{\le -1}E_{1,\lambda} $.
We have the isomorphism for degree
$ \le -1 $ if the last assumption of (0.2) is satisfied, i.e.
if for
$ |I| < n $ we have
$ H^{j}(Y_{I} \cap L_{x},\bC) = 0 $ except for
$ j = 0 $ or
$ 2n + 1 - 2|I| $.

Let
$ K(\bC; v_{1},\dots, v_{m}) $ be the Koszul complex for
$ v_{i} = id : \bC \to \bC \, (1 \le i \le m) $.
Then
$$
\sigma_{\le n-1}K(\bC; v_{1},\dots, v_{m}) =
\sigma_{\le n-1}K_{\lambda},
\leqno(4.3.4)
$$
and
$ \sigma_{\le n}K(\bC; v_{1},\dots, v_{m}) $ is a
direct factor of
$ \sigma_{\le n}K_{\lambda} $, because
$ Y_{I} $ may be reducible if
$ |I| = n $.
So we may replace
$ K_{\lambda} $ with the Koszul complex as long as we
calculate the cohomology of degree
$ \le n - 1 $.
Since this Koszul complex is acyclic and the rank of its
$ j $-th component is
$ \binom mj $, the rank of the nonzero cohomology group of
$ \sigma_{\ge j}K_{\lambda} $ (i.e. the image of the
differential
$ d^{j-1}) $ is
$ \binom {m-1}{j-1} $ for
$ j \le n - 1 $ by the binomial relation.
So the assertion for
$ \lambda \ne 1 $ follows from Theorem (0.1),
where the shift of the index
$ j $ comes from the fact that the complex
$ K_{\lambda} $ is indexed by
$ |I| $ instead of
$ k = |I|-1 $.

For
$ \lambda = 1 $, we use a (generalized) weight spectral
sequence similar to (2.4):
$$
E_{1}^{-k,j+k} =
\mopls_{|I|=k} H^{j-k}(L_{X,x}\cap Y_{I},\bQ)(-k)
\Rightarrow H^{j}(L_{X,x}\setminus Y,\bQ).
$$
This is induced by the weight filtration
$ W $ on
$ (\boR j'_{*}\bQ_{X\setminus Y})|_{Y\setminus \{x\}} $
(see [6]) such that
$$
\Gr_{k}^{W}(\boR j'_{*}\bQ_{X\setminus Y})|_{Y\setminus \{x\}}
= \mopls_{|I|=k} \bQ_{Y_{I}\setminus \{x\}}(-k)[-k],
$$
where
$ j' : X \setminus Y \to X $ is as in (1.7).
So
$ H^{j}(L_{X,x}\setminus Y,\bQ) $ is of type
$ (j,j) $ for
$ j \le n - 1 $
because
$ H^{j}(L_{X,x}\cap Y_{I},\bQ) = 0 $ for
$ j \le n - 1 - |I| $ by (4.1.2) and (1.6).

Since
$ H^{j}(F_{x},\bQ)_{1} $ has weights
$ \le 2j $ and
$ N $ is a morphism of type
$ (-1,-1) $, this assertion implies that
$ N = 0 $ on
$ H^{j}(F_{x},\bQ)_{1} $ for
$ j \le n - 1 $ using the Wang sequence (1.7) and
considering
$ \Ker\, N $.
The assertion on the rank then follows using the Wang sequence
and the binomial relation.
This completes the proof of Theorem (0.2).

\medskip\noindent
{\bf 4.4.~Remark.}
In [9], the case
$ a_{i} = 1 $ for any
$ i $ was treated.
The arguments there (e.g. Th. 3.1) imply also the assertion
on the rank in (0.2) in this case (see also [5], [15] for
the case of a generic central arrangement), and Th. 5.1
corresponds to the vanishing results in (0.2).
In Cor. 4.1, it is proved that the monodromy is trivial for
$ j \le n - 1 $ in this case.

\medskip\noindent
{\bf 4.5.~Proof of Theorem (0.3).}
Let
$ \cF_{\lambda} = \psi_{f,\lambda}\bC_{X}|_{U} $, and let
$ j_{U} : U \to E $ denote the inclusion morphism.
By (1.5.2) we have canonical isomorphisms (compatible with
$ T $)
$$
H^{i}(F_{x},\bC)_{\lambda} =
H^{i}(E,\psi_{f,\lambda}\bC_{X}|_{E}) =
H^{i}(U,\cF_{\lambda}).
\leqno(4.5.1)
$$
By (1.4),
$ \cF_{\lambda} $ is a local system of rank
$ 1 $ if
$ \lambda^{e} = 1 $, and
$ \cF_{\lambda} = 0 $ otherwise.
So the action of the monodromy
$ T $ on
$ \cF_{\lambda} $ and
$ H^{i}(F_{x},\bC)_{\lambda} $ is the multiplication by
$ \lambda $ (i.e. semisimple).
The monodromy of
$ \cF_{\lambda} $ around
$ Y_{j} $ is given by the multiplication by
$ \lambda^{-a_{j}} $.
By (4.5.1) we get
$$
\chi_{\lambda}(F_{x}) = \chi (U)\,\,\, \text{if
$ \lambda^{e} = 1 $, and
$ 0 $ otherwise}.
\leqno(4.5.2)
$$
Since we assume that the
$ Y_{j} $ are principal, we have
$$
H^{j}(F_{x},\bC)_{\lambda} = 0\,\,\, \text{for
$ j \ne n $, if
$ \lambda^{a_{i}} \ne 1 $ for some
$ i $,}
\leqno(4.5.3)
$$
using the weak Lefschetz theorem, because
$ E \setminus Y'_{i} $ is affine where
$ Y'_{i} $ is the proper transform of
$ Y_{i} $.
Indeed the last assertion can be reduced to the case
$ X $ smooth, replacing
$ X $ with an ambient smooth space, because
$ Y_{i} $ is principal.
So we get
$$
b_{\lambda}^{n}(F_{x}) = (-1)^{n}\chi (U)\,\,\,
\text{if
$ \lambda^{e} = 1 $ and
$ \lambda^{d} \ne 1 $.}
\leqno(4.5.4)
$$

If
$ \lambda^{d} = 1 $, then it is known that
$ \cF_{\lambda} $ is a constant sheaf on
$ U $.
(Indeed,
$ \mopls_{\lambda} \cF_{\lambda} $ is the direct image of
a constant sheaf on a finite covering of
$ U $ which is ramified over
$ E \cap Y' $,
see [20], etc.)
Let
$ D_{j} := E \cap Y'_{i} $,
and
$ D^{(k)} $ be the disjoint union of
$ D_{I} := \mcap_{j\in I}D_{j} $ for
$ |I| = k $ where
$ D_{\emptyset} = E $.
Then the cohomology of
$ U $ is calculated by using the weight spectral sequence [6]
$$
E_{1}^{-k,j+k} = H^{j-k}(D^{(k)},\bQ(-k)) \Rightarrow
H^{j}(U,\bQ).
\leqno(4.5.5)
$$

By assumption the constant sheaf
$ \bQ_{X}[n+1] $ is a perverse sheaf, and hence so are
$ \bQ_{Y_{I}}[n+1-|I|] $ for any
$ I $,
where
$ Y_{I} = \mcap_{j\in I} Y_{j} $,
see (1.6).
On the other hand, it is known that, if there is a blow-up
$ \pi : X' \to X $ with a point center such that
$ X' $ and the exceptional divisor
$ E $ are smooth, then the primitive cohomology of
$ E $ is isomorphic to the stalk of the intersection
cohomology
$ \IC_{X}\bQ $ of
$ X $ at
$ x $.
(Indeed, by the decomposition theorem [2],
$ \boR\pi_{*}\bQ_{X'} = \IC_{X}\bQ[-n-1] \oplus
M^{\ssbull} $ with
$ \supp M^{\ssbull} = \{x\} $, and
$ M $ is symmetric with center
$ n+1 $, i.e.
$ \dim M^{n+1-j} = \dim M^{n+1+j} $ by the relative
hard Lefschetz theorem for
$ \pi $.
On the other hand,
$ H^{\ssbull}(E,\bQ) = (\IC_{X}\bQ[-n-1])_{x} \oplus
M^{\ssbull} $, and it is symmetric with center
$ n $ by the classical hard Lefschetz theorem.
Then the assertion follows from the Lefschetz decomposition
because
$ H^{j}(E,\bQ) = M^{j} $ for
$ j > n $.)
So the
$ j $-th primitive cohomology of the exceptional divisor
vanishes for
$ 0 < j < \dim X - 1 $,
using an exact sequence as in (4.1.3).
Similar assertions hold also for any
$ Y_{I} $.

For
$ 0 \le j < n $,
the above arguments imply that
$$
E_{1}^{-k,j+k} = \mopls_{|I|=k} \bQ
(\hbox{$ \frac {-j-k}2 $}),
\leqno(4.5.6)
$$
if
$ j + k $ is even and
$ 0 \le k \le j $,
and it is zero otherwise.
For
$ j = n $,
$ E_{1}^{-k,n+k} $ contains
$ \mopls_{|I|=k} \bQ({-n-k \over 2}) $ if
$ j + k $ is even and
$ 0 \le k \le n $.
Furthermore, the differential is given by the co-Cech Gysin
morphism.
Thus the cohomology of the
$ E_{1} $-complex of (4.5.5) for
$ j < n $ is calculated by that of
$$
\mopls_{j\ge 0}\sigma_{\ge m-j}
K(\bC; v_{1},\dots,v_{m})[m-2j](-j),
$$
where the Koszul complex
$ K(\bC; v_{1},\dots,v_{m}) $ is as in (4.3).
So the assertion on the Hodge type and the rank in (0.2)
holds for
$ j \le n - 1 $ in this case.

Combined with (4.5.4), this implies
$$
b_{\lambda}^{n}(F_{x}) = (-1)^{n}\chi (U) +
\hbox{$ \binom {m-2}{n-1} $} \,\,\,
\text{if
$ \lambda^{e} = 1 $ and
$ \lambda^{d} = 1 $,}
\leqno(4.5.7)
$$
because
$ \binom {m-2}{n-1} = \msum_{0\le k\le n-1}
(-1)^{k}\binom {m-1}{n-1-k} $, see [15], lemma 2.5.
This completes the proof of Theorem (0.3).

\medskip\noindent
{\bf 4.6.~Remark.}
With the assumption of (4.5), assume further
$ X $ smooth.
Then it is known that
$ \chi (U) $ is explicitly calculated by using
$ d_{j} $.
Indeed, we have by (4.5.5)
$$
\chi (U) = \mopls_{|I|\le n} (-1)^{|I|}\chi (D_{I}).
$$
Furthermore, by the theory of Chern classes (see e.g. [11]),
the topological
Euler characteristic
$ \chi (D_{I}) $ is the coefficient of
$ T^{n} $ in
$$
(1+T)^{n+1}\mprod_{j\in I} (d_{j}T/(1+d_{j}T)) \in
\bQ[[T]]/(T^{n+1}),
$$
because the
$ k $-th Chern class of the tangent bundle of
$ D_{I} $ gives the topological Euler characteristic for
$ k = \dim D_{I} \,(= n - |I|) $,
and the restriction of a cycle on
$ \bP^{n} $ to
$ D_{I} $ is essentially same as the intersection with
$ D_{I} $.
Here the truncated formal power series ring is identified
with the cohomology ring of
$ \bP^{n} $ so that
$ (1+T)^{n+1} $ is the total Chern class of the tangent
bundle of
$ \bP^{n} $,
and
$ 1+d_{j}T $ is that of the normal bundle of
$ D_{j} $.

Since
$ 1 - d_{j}T/(1+d_{j}T) = (1+d_{j}T)^{-1} $, we see that
$ \chi (U) $ is the coefficient of
$ T^{n} $ in
$$
(1+T)^{n+1}\mprod_{1\le j\le m} (1+d_{j}T)^{-1} \in
\bQ[[T]]/(T^{n+1}).
$$
For
$ m = 1 $ and
$ a_{1} = 1 $, this is compatible with a well-known formula
for the Milnor number of a homogeneous hypersurface isolated
singularity (using Theorem (0.3)), i.e.
$$
1 - d_{1}\chi(U) = (1-d_{1})^{n+1}.
$$
In the case of a generic central arrangement (i.e.
$ d_{j} = 1 ), $ the above assertion implies
$$
\chi(U) = (-1)^{n} \hbox{$ \binom {m-2}n $}.
\leqno(4.6.1)
$$
This is compatible with the formula in [5], [15]
using Theorem (0.3).

In general, we can verify that the coefficient of
$ T^{k} $ in
$$
\mprod_{1\le j\le m} (1+d_{j}T)^{-1} \in
\bQ[[T]]/(T^{n+1})
$$
is a polynomial in
$ d_{1},\dots , d_{m} $,
which is equal to
$$
(-1)^{k} \msum_{1\le i\le m} ({d}_{i}^{k+m-1}/
\mprod_{p\ne i} (d_{i}-d_{p}))
$$
in the fraction field
$ \bQ(d_{1},\dots , d_{m}) $.
This follows by induction on
$ m $, using
$$
\msum_{0\le j\le k} {d}_{m-1}^{j}{d}_{m}^{k-j} =
({d}_{m-1}^{k+1} - {d}_{m}^{k+1})/(d_{m-1} - d_{m}).
$$
Furthermore, the above polynomial vanishes for
$ 1-m \le k < 0 $,
because it is a polynomial, and has negative degree.
So we see that
$ \chi (U) $ for
$ m > 1 $ is a polynomial in
$ d_{1},\dots , d_{m} $,
which is equal to
$$
\msum_{1\le i\le m} (- {d}_{i}^{m-2}(1-d_{i})^{n+1}/
\mprod_{p\ne i} (d_{i}-d_{p}))
\leqno(4.6.2)
$$
in the fraction field.
This gives an explicit formula if the
$ d_{i} $ are different from each other.
In general we have to take a limit
(or make some calculation in the fraction field).

\bigskip\bigskip
\centerline{\bf Appendix}

\bigskip\noindent
We give an example such that the monodromy at degree
$ n - 1 $ is not semisimple at the origin, but is
semisimple at the other points.
This shows that Theorem (0.4) is optimal, and that the
extension class between the graded-pieces of the filtration
associated to the Leray spectral sequence (2.2) for the
nearby cycles is nontrivial as
$ \bC[N] $-modules.

\medskip\noindent
{\bf A.1.~Embedded resolution of singularities.}
We first explain how to get an embedded resolution of
a function of type
$ f = f_{d} + f_{d+1} $ on the affine cone
$ X $ of a smooth projective variety
$ E $ with a very ample line bundle
$ L $ defining the embeddings
$ E \to \bP^{r-1} $ and
$ X \to \bC^{r} $ where
$ r = \dim \Gamma(E,L) $.
Here
$ f_{j} $ is an element of the
$ j $-th symmetric power of
$ \Gamma(E,L) $, which is identified with a polynomial
of degree
$ j $ in
$ r $ variables, and defines a function on the affine
cone
$ X $.
We assume that
$ f_{d}^{-1}(0) \setminus f_{d+1}^{-1}(0) $ defines a
divisor with simple normal crossings on
$ E \setminus f_{d+1}^{-1}(0) $ (where
$ f_{j} $ is viewed as a section of
$ L^{\otimes j} $).

Let
$ X^{\vee} $ be the total space of the dual of the line
bundle
$ L $ with the projection
$ \rho : X^{\vee} \to E $.
It is the blow-up of the affine cone
$ X $ at the origin, and the exceptional divisor is
identified with
$ E $.
Let
$$
D_{0} = f_{d}^{-1}(0),\quad
D_{\infty} = f_{d+1}^{-1}(0)
$$
as (not necessarily reduced) divisors on
$ E $.
Let
$ Y^{\vee} $ be the proper transform of
$ f^{-1}(0) $ in
$ X^{\vee} $.
Let
$ D_{0,\infty} $ be the greatest common divisor of
$ D_{0} $ and
$ D_{\infty} $, and put
$$
D_{0}^{\red} = D_{0} - D_{0,\infty},\quad
D_{\infty}^{\red} = D_{\infty} - D_{0,\infty} .
$$
Then we have a canonical decomposition
$$
Y^{\vee} = Y_{\hor} + Y_{\ver},
$$
where
$ Y_{\ver} = \rho^{*}D_{0,\infty} $ and
$ Y_{\hor} $ corresponds to a rational section
$ \sigma $ of the line bundle such that
$$
\div \sigma = D_{0}^{\red} - D_{\infty}^{\red}.
$$

Assume there is an embedded resolution
$ \pi : E' \to E $ of
$ D_{0} $ such that
$ D_{0} \cup D_{\infty} $ is a divisor with normal
crossings on a neighborhood of
$ D_{0} $.
(This is satisfied in the case
$ n = 2 $.)
Let
$ \pi : X' \to X^{\vee} $ be the base change of
$ \pi : E' \to E $ by
$ \rho $.
We can similarly define
$ D'_{0,\infty} $,
$ D_{0}^{\prime\,\red} $,
$ D_{\infty}^{\prime\,\red} $,
$ Y'_{\hor}, Y'_{\ver} $ for
$ D'_{0} = \pi^{*}D_{0} $,
$ D'_{\infty} = \pi^{*}D_{\infty} $,
$ Y' = \pi^{*}Y^{\vee} $ so that
$$
\div \pi^{*}\sigma =
D_{0}^{\prime\,\red} - D_{\infty}^{\prime\,\red}.
$$

Blowing up further if necessary, we may assume
$$
D_{0}^{\prime\,\red} \cap D_{\infty}^{\prime\,\red} =
\emptyset.
\leqno\text{\rm (A.1.1)}
$$
Then we get an embedded resolution of
$ f^{-1}(0) $ by iterating blow-ups of
$ X' $ along the irreducible components of
$ D_{0}^{\prime\,\red} $.
Indeed,
$ Y'_{\hor} $ may be locally defined by
$ s = \prod_{i}x_{i}^{m_{i}} $ with
$ x_{1},\dots, x_{n} $ local coordinates of
$ E' $ and
$ s $ a local coordinate of the line bundle so that
the blow-up along
$ \{x_{i} = s = 0\} $ corresponds to the substitution of
$ s $ by
$ sx_{i} $ where
$ m_{i} $ decreases by
$ 1 $.

For simplicity, assume
$ n = 2 $,
$ D_{0} $ is a reduced divisor with simple normal
crossings, and intersects
$ D_{\infty} $ at smooth points of
$ D_{0} $.
Since the embedded resolution can be obtained by iterating
blowing-ups with point centers, we can verify that
$ D_{0}^{\prime\,\red} $ may be assumed to be isomorphic to
$ D_{0} $,
and does not intersect
$ D_{\infty}^{\prime\,\red} $ (calculating the multiplicities
of the exceptional divisors).
If furthermore
$ D_{0} $ is smooth, then the exceptional divisor of the
blow-up along
$ D_{0}^{\prime\,\red} $ is a trivial
$ \bP^{1} $-bundle over
$ D_{0}^{\prime\,\red} $,
because the proper transform of
$ Y^{\vee} $ gives a trivialization.

For example, if
$ D_{0} $ (resp.
$ D_{\infty}) $ is defined locally by
$ y = 0 $ (resp.
$ x = 0) $ with multiplicity
$ 1 $ (resp.
$ m) $, then the resolution is obtained by iterating
$ m $ times blow-ups along a point of the proper transform of
$ D_{0} $.
Let
$ C_{j} $ denote the proper transform of the exceptional
divisor of the
$ j $-th blow-up for
$ 1 \le j \le m $.
Then
$$
\pi^{*}D_{0} =
\msum_{1 \le j \le m} jC_{j} + D_{0}^{\prime\,\red},
\quad
\pi^{*}D_{\infty} =
\msum_{1 \le j \le m} mC_{j} + D''_{\infty},
$$
where
$ D''_{\infty} $ is the proper transform of
$ D_{\infty} $ (with multiplicity
$ m) $, and
$$
D'_{0,\infty} = \msum_{1 \le j \le m} jC_{j},
\quad
D_{\infty}^{\prime\,\red} =
\msum_{1 \le j \le m-1} (m-j)C_{j} + D''_{\infty}.
$$

\medskip\noindent
{\bf A.2.~Conditions for non semisimplicity.}
With the notation and the assumptions of (1.5), assume
$ n = 2 $.
We consider the conditions for
the non semisimplicity of the monodromy on
$ H^{1}(F_{x},\bQ)_{\lambda} $.
Define
$$
\aligned
J(\lambda;a,b)
&= \{I \subset J(\lambda) : |I| - 1 = a, s(I) \ge b\},
\\
J_{0}(\lambda;a,b)
&= \{I \in J(\lambda;a,b) : Y_{I} \cap Y_{j} = \emptyset
\,\,\,\text{for}\,\, j \notin J(\lambda) \}.
\endaligned
$$
Let
$ u $ be an element of
$ E_{1}^{-1,2} $ in (1.5.3).
It may be viewed as an element of
$$
\mopls_{I\in J_{0}(\lambda;1,0)} H^{0}(Y_{I},\bQ),
$$
because
$ H^{0}(U_{I},\cF_{\lambda,I}) $ vanishes for
$ I \in J(\lambda;1,0) \setminus J_{0}(\lambda;1,0) $.
Here the Tate twist
$ (-1) $ is trivialized by choosing
$ \sqrt{-1} $.

The first condition on
$ u $ is that it is annihilated by the differential
$ d_{1} $ of the spectral sequence, i.e.
its images in
$$
\mopls_{I\in J_{0}(\lambda;2,1)}
H^{0}(Y_{I},\bQ),\quad
\mopls_{I\in J(\lambda;0,0)}
H^{2}(U_{I},\cF_{\lambda,I})(1)
$$
vanish.
This condition is necessary to assure that it defines an
element of
$ \Gr_{2}^{W}H^{1}(F_{x},\bQ)_{\lambda} $.

The second condition is that its image in
$ \mopls_{I\in J_{0}(\lambda;1,1)} H^{0}(Y_{I},\bQ) $
does not belong to the image of
$ \mopls_{I\in J_{0}(\lambda;0,0)} H^{0}(Y_{I},\bQ) $.
This condition is necessary to assure that its image by
$ N $ does not vanish in
$ \Gr_{0}^{W}H^{1}(F_{x},\bQ)_{\lambda} $.

\medskip\noindent
{\bf A.3.~Example.}
$ $ Let
$ X = \{xw - yz = 0\} \subset \bC^{4} $,
and
$ f = (y^{2} - x^{4})(x^{2} - y^{4}) $ where
$ n = 2 $.
Then
$ E = \bP^{1} \times \bP^{1} $ with coordinates
$ (u_{0},u_{1};v_{0}, v_{1}) $ such that
$$
x = u_{0}v_{0},\quad
y = u_{0}v_{1},\quad
z = u_{1}v_{0},\quad
w = u_{1}v_{1}.
$$
We apply the arguments in (A.1) to
$$
g = x \pm y^{2},\quad
h = y \pm x^{2},
$$
where
$ g_{1} = x, g_{2} = y^{2} $, etc., and
$ L $ is induced by
$ \cO(1) $ on
$ \bP^{3} $.
Let
$ (u,v) $ be the affine coordinates on
$ \{u_{1}v_{1} \ne 0\} \subset E $ such that
$ u = u_{0}/u_{1}, v = v_{0}/v_{1} $.
Let
$ s $ be the coordinate of the line bundle over
$ \{w \ne 0\} $, which is induced by
$ w $.
Then
$ Y^{\vee} $ near
$ (0,1;0,1) $ is locally defined by
$$
u^{4}(v - su)(v + su) = 0,
\leqno\text{\rm (A.3.1)}
$$
because
$ g = w(x/w) \pm w^{2}(y/w)^{2} $ (and similarly for
$ h $), where
$ w $ is actually
$ s $.
We have a similar assertion on a neighborhood of
$ (0,1;1,0) $.
So
$ Y^{\vee} $ has four reduced components (defined by
$ v \pm su = 0 $, etc.) and one multiple component
(defined by
$ u^{4} = 0 $).

Let
$ Z_{1}, Z_{2} $ be the divisors defined by
$ v_{0} $ and
$ v_{1} $ respectively.
Then
$ D_{0}^{\red} $ in (A.1) for
$ g $ (resp.
$ h $) is
$ Z_{1} $ (resp.
$ Z_{2} $), and
$ D_{0}^{\red} \cap D_{\infty}^{\red} $ consists of
$ (0,1;0,1) $ (resp.
$ (0,1;1,0) $).
Let
$ \pi : E' \to E $ be the blow-up along these two
points with exceptional divisors
$ C_{1}, C_{2} $.
This gives a resolution satisfying (A.1.1) by the last
argument of (A.1) where
$ m = 1 $.
Let
$ Z'_{1}, Z'_{2} $ be the proper transforms of
$ Z_{1}, Z_{2} $ so that
$$
\pi^{*}Z_{1} = Z'_{1} + C_{1},\quad
\pi^{*}Z_{2} = Z'_{2} + C_{2}.
\leqno\text{\rm (A.3.2)}
$$

Let
$ \pi : X' \to X^{\vee} $ be the base change of
$ \pi : E' \to E $ by
$ \rho $.
Let
$ X'' \to X' $ be the blow-up along
$ Z'_{1} $ and
$ Z'_{2} $ with exceptional divisors
$ E_{1}, E_{2} $.
This gives an embedded resolution of
$ f^{-1}(0) $.
We see that
$ E_{1} $ is a trivial
$ \bP^{1} $-bundle over
$ Z'_{1} $, and the intersection of
$ E_{1} $ with the proper transform of
$ f^{-1}(0) $ consists of two connected components
(corresponding to
$ v - su = 0 $ and
$ v + su = 0) $ and these are both isomorphic to
$ Z'_{1} $ by the projection (and similarly for
$ E_{2}, Z'_{2} $).
Let
$ E_{0} $ be the proper transform of the zero section
$ E' $ by
$ X'' \to X' $.
For
$ i = 1, 2 $, the proper transform of
$ \rho^{-1}(C_{i}) $ will be denoted by
$ E_{i+2} $.
Let
$ C'_{i} $ be the proper transform of
$ C_{i} $, which is equal to
$ E_{0} \cap E_{i+2} $.
We will identify
$ Z'_{i} $ with
$ E_{0} \cap E_{i} $ for
$ i = 1, 2 $.
Note that the inverse image of the origin is
$ \mcup_{0\le i \le 2}E_{i} $.

Using this resolution together with the conditions in
(A.2), we can show that the action of
$ N $ on
$ H^{1}(F_{0},\bQ)_{\lambda} $ is not semisimple where
$ \lambda = -1 $.
We see that the multiplicities of the irreducible
components are even except for the proper transforms of
the four reduced components of
$ Y^{\vee} $.
We have to find an appropriate element
$ u $ as in (A.2).
We define
$ u $ by
$$
1 \in H^{0}(Z'_{1},\bQ),\quad
1 \in H^{0}(C'_{1},\bQ),\quad
-1 \in H^{0}(Z'_{2},\bQ),\quad
-1 \in H^{0}(C'_{2},\bQ).
$$
Here we use the natural order of the exceptional divisors
$ E_{i} $ for
$ 0 \le i \le 4 $ to define these elements, because Cech and
co-Cech complexes are involved.
We can verify that the two conditions in (A.2) are satisfied
by using (A.3.2), etc.
Note that, if
$ Y_{I} $ is
$ E_{i} $ with
$ i = 1 $ or
$ 2 $, then
$$
H^{j}(U_{I},\cF_{\lambda,I}) = 0\quad
\text{for any}\,\, j,
$$
because
$ U_{I} $ is the product of
$ Z'_{i} $ with
$ \bP^{1} $ minus two points, and the monodromy of
$ \cF_{\lambda,I} $ around the two points are
$ -1 $ (here we use the Leray spectral sequence for
the projection to
$ Z'_{i} $).
We can also verify that the Milnor monodromy is semisimple
outside the origin, using (A.3.1) and (1.5).

\medskip\noindent
{\bf A.4.~Remark.}
For the moment, we do not know any example as above with
$ X $ smooth.

\bigskip\bigskip
\centerline{{\bf References}}

\bigskip

\item{[1]}
D. Barlet, Interaction de strates cons\'ecutives pour les
cycles \'evanescents, Ann. Sci. Ecole Norm. Sup. (4) 24 (1991),
401--505.

\item{[2]}
A. Beilinson, J. Bernstein and P. Deligne, Faisceaux pervers,
Ast\'erisque, vol. 100, Soc. Math. France, Paris, 1982.

\item{[3]}
J.L. Brylinski, Transformations canoniques, dualit\'e
projective, th\'eorie de Lefschetz, transformations de
Fourier et sommes trigonom\'etriques, in G\'eom\'etrie
et analyse microlocales, Ast\'erisque 140--141 (1986), 3--134.

\item{[4]}
N. Budur and M. Saito, Multiplier ideals,
$ V $-filtration, and spectrum, preprint (math.AG/0305118).

\item{[5]}
D. Cohen and A. Suciu, On Milnor fibrations of arrangements,
J. London Math. Soc. 51 (1995), 105--119.

\item{[6]}
P. Deligne, Th\'eorie de Hodge I, Actes Congr\`es
Intern. Math., 1970, vol. 1, 425-430; II, Publ. Math. IHES,
40 (1971), 5--57; III, ibid., 44 (1974), 5--77.

\item{[7]}
\SameAuthor, Le formalisme des cycles \'evanescents, in
SGA7 XIII and XIV, Lect. Notes in Math. vol. 340, Springer,
Berlin, 1973, pp. 82--115 and 116--164.

\item{[8]}
A. Dimca, Sheaves in Topology, Universitext, Springer, Berlin
2004.

\item{[9]}
A. Dimca and A. Libgober,
Local topology of reducible divisors, preprint (math.AG/
0303215).

\item{[10]}
A. Dimca and M. Saito, Monodromy at infinity and the weights of
cohomology, Compos. Math. 138 (2003), 55--71.

\item{[11]}
W. Fulton, Intersection theory, Springer, Berlin, 1984.

\item{[12]}
M. Kashiwara, Vanishing cycle sheaves and holonomic systems
of differential equations, Algebraic geometry (Tokyo/Kyoto,
1982), Lect. Notes in Math. 1016, Springer, Berlin,
1983, pp. 134-142.

\item{[13]}
B. Malgrange, Polyn\^ome de Bernstein-Sato et cohomologie
\'evanescente, Analysis and topology on singular spaces, II,
III (Luminy, 1981), Ast\'erisque 101--102 (1983), 243--267.

\item{[14]}
V. Navarro Aznar, Sur la th\'eorie de Hodge-Deligne, Inv.
Math. 90 (1987), 11--76.

\item{[15]}
P. Orlik and R. Randell, The Milnor fiber of a generic arrangement,
Ark. Mat. 31 (1993), 71--81.

\item{[16]}
M. Saito, Modules de Hodge polarisables, Publ. RIMS, Kyoto
Univ. 24 (1988), 849--995.

\item{[17]}
\SameAuthor, Mixed Hodge Modules, Publ. RIMS, Kyoto Univ. 26
(1990), 221--333.

\item{[18]}
D. Siersma, Variation mappings on singularities with a
1-dimensional critical locus, Topology 30 (1991), 445--469.

\item{[19]}
\SameAuthor, The vanishing topology of non isolated singularities,
in New Developments in Singularity Theory (D. Siersma et al.,
eds.), Kluwer Acad. Publishers, 2001, pp. 447--472.

\item{[20]}
J.H.M. Steenbrink, Mixed Hodge structure on the vanishing
cohomology, in Real and Complex Singularities (Proc. Nordic
Summer School, Oslo, 1976) Alphen a/d Rijn, Sijthoff \&
Noordhoff 1977, pp. 525--563.

\bigskip
Alexandru Dimca

Laboratoire J.A. Dieudonn\'e, UMR du CNRS 6621
                 
Math\'ematiques
                 
Universit\'e de Nice-Sophia-Antipolis
                 
Parc Valrose,
          
06108 Nice Cedex 02, FRANCE. 

e-mail: dimca\@math.unice.fr

\medskip
Morihiko Saito

RIMS Kyoto University, Kyoto 606--8502 JAPAN

e-mail: msaito\@kurims.kyoto-u.ac.jp

\bigskip
\vers

\bye